\documentclass[letterpaper, 10 pt, conference]{ieeeconf}  % Comment this line out if you need a4paper

\IEEEoverridecommandlockouts                              % This command is only needed if 
\usepackage{graphics} % for pdf, bitmapped graphics files
\usepackage{graphicx}
\usepackage{epsfig} % for postscript graphics files
\usepackage{times} % assumes new font selection scheme installed
\usepackage{amsmath,amssymb} % assumes amsmath package installed
\newcommand{\R}{\mathbb{R}}

\DeclareMathOperator{\argmin}{argmin}
\usepackage[linesnumbered,ruled]{algorithm2e}
\usepackage[justification=centering]{caption}
\usepackage{flushend}
\allowdisplaybreaks

\title{\LARGE \bf Multi Agent Pathfinding for Noise Restricted Hybrid Fuel Unmanned Aerial Vehicles}

\author{Drew Scott$^{1}$, Satyanarayana G. Manyam$^{2}$, David W. Casbeer$^{3}$, Manish Kumar$^{4}$, Isaac E. Weintraub$^{5}$% <-this % stops a space
\thanks{$^{1}$Department of Mechanical and Materials Engineering, University of Cincinnati 
        {\tt\small scott2dd@mail.uc.edu}}%
\thanks{$^{2}$Research Scientist, Infoscitex corp., a DCS company, Dayton OH, 45431
        {\tt\small msngupta@gmail.com}}%
\thanks{$^{3}$Technical Area Lead, Cooperative \& Intelligent Control, Control Science
Center, Air Force Research Laboratory, WPAFB, OH, 45433
        {\tt\small david.casbeer@us.af.mil}}%
\thanks{$^{4}$Professor in Department of Mechanical and Materials Engineering, University of Cincinnati {\tt\small kumarmu@ucmail.uc.edu}}
\thanks{$^{5}$Electronics Engineer, Control Science Center, Air Force Research Laboratory, WPAFB, OH, 45433 {\tt\small isaac.weintraub.1@us.af.mil}}%
\thanks{This material is in part based on research sponsored by the Ohio Department of Higher Education and the Southwestern Council for Higher Education under Ohio House Bill 49 of the 132nd General Assembly.  The U.S. Government is authorized to reproduce and distribute reprints for Governmental purposes notwithstanding any copyright notation thereon.  The views and conclusions contained herein are those of the authors and should not be interpreted as necessarily representing the official policies or endorsements, either expressed or implied, of Southwestern Council for Higher Education, the Ohio Department of Higher Education or the U.S. Government.}
\thanks{Distribution Statement A. Approved for public release, distribution unlimited. Case Number: AFRL-2023-4807.}
}

\SetAlgoSkip{medskip}
\begin{document}
\maketitle
\thispagestyle{empty}
\pagestyle{empty}
%%%%%%%%%%%%%%%%%%%%%%%%%%%%%%%%%%%%%%%%%%%%%%%%%%%%%%%%%%%%%%%%%%%%%%%%%%%%%%%%%%%%%%%%%%%%%%%%

\begin{abstract}
Multi Agent Path Finding (MAPF) seeks the optimal set of paths for multiple agents from respective start to goal locations such that no paths conflict.  We address the MAPF problem for a fleet of hybrid-fuel unmanned aerial vehicles which are subject to location-dependent noise restrictions.  We solve this problem by searching a constraint tree for which the subproblem at each node is a set of shortest path problems subject to the noise and fuel constraints and conflict zone avoidance.  A labeling algorithm is presented to solve this subproblem, including the conflict zones which are treated as dynamic obstacles.  We present the experimental results of the algorithms for various graph sizes and number of agents.
\end{abstract}

\section{INTRODUCTION}\label{sec:intro}
Hybrid-fuel Unmanned Aerial Vehicles (UAVs) are those in which multiple fuel sources are used in combination to provide energy storage and power for the vehicle.  Currently, these are most often utilized in applications which require extreme flight range in low powered systems.  However, hybrid-fuel platforms will likely see an increase in popularity in the future as UAVs in general become more widely used and thus the scenarios under which hybrid-fuel UAVs are advantageous also grow.  While the primary advantage of hybridization is increased flight-range, if power modality is able to be switched throughout the mission, power management can be employed to exploit the advantages of each energy source.  For example, a combustion or jet engine has high energy-endurance, but is extremely noisy, while a battery-pack provides low energy-density but allows quieter operation.  In the future, we envision noise-production by low-altitude UAVs being a major concern for both operators and the general public.  A survey of the public \cite{yedavalli2019assessment} found that a primary concern with increased small UAV usage was the noise production - specifically the general volume, duration, and time of day.  These concerns are likely to be manifest as UAVs become more widely used for a variety of tasks.  Thus, we envision restrictions on the noise produced by UAVs in certain airspaces.  In the case of urban path-planning, the noise levels become increasingly important as the average number of UAVs overhead increases.  These restrictions may limit the use of engine-powered flights over certain sections of an airspace. 

We consider here a fleet of hybrid-fuel UAVs for which the propeller is powered by a battery-pack, with a combustion engine gen-set onboard to recharge the battery.  Such a platform is defined in \cite{TOWNSEND2020e05285} as a series-hybrid platform.  In the case of these vehicles, the engine may be turned off to allow travel through noise-restricted areas for which the engine is too noisy, but flying is a battery-only mode is quiet enough for travel through such a zone.  The high energy density of traditional fuel allows longer endurance over battery storage alone.  Thus, these series-hybrid UAVs achieve extended flight range via the combustion engine while being capable of flying through noise-restricted zones that are too restrictive for engine-powered flight.  This scenario couples the path planning and power management, such that a path and power plan must be found in tandem. 

The single-agent version of this problem, referred to as the Noise-Restricted Hybrid-Fuel Shortest Path Problem (NRHFSPP), is studied in \cite{scott2024power}. We extend that single agent problem to the Multi-Agent Path-Finding (MAPF) problem with the series-hybrid UAV fleet and noise restrictions, which we refer to as Noise-Restricted Hybrid-Fuel MAPF (NRHF-MAPF).  The primary motivation of the noise restrictions is widespread use of UAVs, particularly in congested areas.  In a congested airspace, there is the immediate concern of cooperative routing and collisions with other aircraft.  Thus, a natural extension of the single-agent problem is the multi-agent path finding extension. 

The MAPF, which is a fundamental problem in multi-agent planning, involves a set of agents which each must be routed from a start location to a goal such that there are no conflicts or collisions between the paths.  The normal objective is to minimize the total travel cost across the agents, but may also include minimizing maximum cost (across agents), or minimizing time to completion.  A usual technique is to discretize the time and environment space, where the latter is modeled as a graph.  Collisions may be defined in a variety of ways, as outlined in \cite{stern2019multi}.  The two main conflict types are i) vertex conflicts, where two or more agents occupy the same vertex at the same time ii) edge conflicts, where 2 or more agents travel along the same edge at the same time.  These conflict types are the one considered here.  

MAPF was shown to be PSPACE-hard in \cite{hopcroft1984complexity,reif1979complexity} and NP-hard in \cite{surynek2010optimization, yu2013structure}.  A global A*-variant search scales poorly as number of agents grows \cite{standley2010finding}.  A two-level search was given in \cite{sharon2013increasing}, referred to as Increasing Cost Tree Search, which was feasible for problems significantly larger than the limit for global A*. Conflict Based Search (CBS) was presented in \cite{sharon2015conflict}.  It is an optimal algorithm that iteratively searches a Constraint Tree (CT) where conflict constraints are added as are violated within subproblems.  Enhanced Conflict Based Search (ECBS) \cite{barer2014suboptimal} is a heuristic algorithm which searches the same CT as CBS, obtaining near-optimal solutions in a scalable manner.   

As stated above, the NRHF-MAPF is an extension of the authors' prior work on single-agent consideration of series-hybrid UAVs in presence of noise-restricted airspaces \cite{scott2022power, scott2024power, manyam2022path, scott2022development, scott2023OCP}.  A similar single-agent problem is studied by another group in \cite{jadischke2023optimal}.  Extension of the standard MAPF for a fleet of hybrid-fuel UAVs involves adding the energy constraints and restrictions on the generator within the noise-restricted zones.  This results in a novel MAPF variant, where conflict free paths must be found for each agent in the system while also respecting the power and noise constraints.  To solve this, we utilize Conflict Based Search (CBS) and replace the standard MAPF subproblem of shortest path problem (SPP) with our NRHFSPP, adapting the algorithm of \cite{scott2024power} to include conflict-zone avoidance. 

The primary contributions of this paper are: (i) presentation and study of the NRHF-MAPF, a novel MAPF variant, (ii) A Mixed Integer Linear Programming (MILP) formulation of the problem, (iii)  Adaptation of an existing NRHFSPP algorithm to deal with dynamic obstacles to use as a subroutine in CBS, and (iv) numerical testing of the presented solution.  

The paper is organized as follows.  The MILP formulation of the problem is given in \ref{sec:prob}.  The CT and CBS for this problem are discussed in Section \ref{sec:CBS}.  The labeling algorithm used to solve the subproblem is presented in Section \ref{sec:labelalgo}.  Results of numerical testing are given in \ref{sec:results}.  Finally, concluding thoughts and areas of future work are given in Section \ref{sec:conclusions}.

\section{Problem Formulation}\label{sec:prob}
We present the MILP formulation of the noise-restricted hybrid-fuel MAPF problem.  

\subsection{Constant Values}
% switch (i,j) to math modes....
%move E/N to top.... (ij) \in E
$D_{ij} \in \R$ - objective cost of edge (i,j)

$T_{ij} \in \R$ - time to travel along edge (i,j)

$C_{ij} \in \R$ - energy cost (battery discharge) of edge (i,j)

$Z_{ij} \in \R$ - battery recharge by generator along edge (i,j)

$G_{ij} \in \{0,1\}$ - is generator allowed be run along edge (i,j) 

$U$ - set of UAVs

$N$ - set of nodes

$E$ - set of edges 

$S_k \in S$ - starting Node for UAV $k$

$F_k \in S$ - goal Node for UAV $k$

$M$ - large constant value

$Q^k_0$ - initial generator state for UAV $k$

$B^k_0$ - initial battery state for UAV $k$

$T^k_0$ - initial time state for UAV $k$

$B^k_{max}$ - Maximum battery state for UAV $k$

\subsection{Decision Variables}

$x^k_{ij} \in \{0,1\}$ - if edge $(i,j) \in E$ used by UAV $k \in U$

$g^k_{ij} \in \{0,1\}$ - if UAV $k$ runs generator on edge (i,j)

$b^k_i \in \R{}$ - state of battery for UAV $k$ at node $i$

$q^k_i \in \R{}$ - state of generator for UAV $k$ at node $i$

$t^k_i \in \R{}$ - time node $i$ is reached by UAV $k$

\subsection{MILP Formulation}
In the following formulation, the set $\delta^+(i)$ and $\delta^-(i)$ represents the set of outgoing and incoming edges from node $i$, respectively.  That is,  $ \delta^+(i) = \{ (i,j) \in E, \, \forall j \in N\}$, and  $ \delta^-(i) = \{ (j,i) \in E, \, \forall j \in N\}$.
% We define a function $\delta^+(\cdot)$ such that $\delta(i)$ gives the set of all nodes which are neighbors of node $i$.  That is, $(i,j) \in E \quad  \forall j \in \delta(i)$.
\begin{align}
    & J = \min_{x} \sum_{k \in U} \sum_{(i,j) \in E}  D_{ij} x^k_{ij}   \label{MILPobj}\\
    & \sum_{e \in \delta^+(S_k)} x^k_{e} = 1 && \forall k \in U      \label{MILPdeg1}\\
    & \sum_{e \in \delta^-(F_k)} x^k_{e} = 1  && \forall k \in U       \label{MILPdeg2}\\
    & \sum_{e \in \delta^-(i)} x^k_{e} - \sum_{e \in \delta^+(i)} x^k_{e} = 0 \quad   &&\forall i \in N\setminus \{S_k,F_k\}, k \in U   \label{MILPdeg3} \\
    & B^k_{max} \geq b^k_j \geq 0 && \forall j \in N \setminus \{S_k\}, k \in U \label{MILPbatt1}\\
    & b^k_{S_k} = B^k_0 && \forall k \in U \label{MILPbatt2} \\
    & b^k_j \leq b^k_i - C_{ij} \nonumber \\  & \quad + Z_{ij} g^k_{ij} + M(1 - x^k_{ij})   && \forall (i,j) \in E, k \in U \label{MILPbatt3}\\
    & b^k_j \geq b^k_i - C_{ij} \nonumber \\  & \quad + Z_{ij} g^k_{ij} - M(1 - x^k_{ij})   && \forall (i,j) \in E, k \in U \label{MILPbatt4}\\
    & q^k_j \geq 0 && \forall i \in N \setminus \{S_k\},\forall k \in U \label{MILPgen1}\\
    & q^k_{S_k} = Q^k_0 && \forall k \in U  \label{MILPgen2}\\
    & q^k_j \leq q^k_i - Z_{ij}g^k_{ij} \nonumber \\& \quad \quad + M(1- x^k_{ij}) && \forall (i,j) \in E, k \in U \label{MILPgen3}\\
    & g^k_{ij} \leq x^k_{ij} && \forall (i,j) \in E, k \in U \label{MILPmisc1}\\
    & g^k_{ij} \leq G_{ij} && \forall (i,j) \in E, k \in U \label{MILPmisc2}  \\
    & t^k_{S_k} = T^k_0 && \forall k \in U  \label{MILPtime1}\\ 
    & t^k_j \geq t^k_i + T_{ij} + M(1- x^k_{ij}) && \forall (i,j) \in E, k \in U \label{MILPtime2}\\
    &&& \hspace{-126pt} |t^a_i - t^b_i| \geq \epsilon  \qquad \qquad \forall i \in N, a \in U, b \in U, a \neq b \label{MILPvertex}  \\
    &&& \hspace{-126 pt} |t^a_j - t^b_i| \geq \epsilon - \epsilon(2 - x^a_{ij} - x^b_{ji}) \nonumber \\ &&& \hspace{-115pt } \forall (i,j) \in E : (j,i) \in E, a \in U, b \in U, a \neq b \subset U \label{MILPedge} 
\end{align}
where $\epsilon$ is some known threshold for time difference with respect to the conflicts.  Equations \eqref{MILPvertex} and \eqref{MILPedge} can be trivially linearized with auxiliary variables and big-M technique \cite{williams2013model}; but is not shown here for brevity of the MILP presentation.  The objective function is given in \eqref{MILPobj}, which aims to minimize total travel cost across all UAVs. The can easily be replaced for other objectives such as makespan or a minmax cost.  Degree constraints are given in \eqref{MILPdeg1}-\eqref{MILPdeg3}.  Battery constraints and generator constraints are given in \eqref{MILPbatt1}-\eqref{MILPgen3}. Equation \eqref{MILPmisc1} constrains the generator to be run on an edge only if that edge is traveled by the UAV, and \eqref{MILPmisc2} restricts the generator from being run in noise-restricted zones. Time tracking is done in \eqref{MILPtime1}-\eqref{MILPtime2}, and vertex and edge deconflicting is enforced by\eqref{MILPvertex} and \eqref{MILPedge}, respectively.

As stated above, the problem is posed on a graph $(N, E)$ where $N$ is the set of nodes and $E$ is the set of edges connecting them, where an edge connecting node $i$ to $j$ is notated by $e = (i,j)$.

\section{Constraint Tree and Conflict Based Search}\label{sec:CBS}
A recent method in the MAPF literature to deal with the conflict constraints is the searching of a Constraint tree (CT), which is the basis of CBS and ECBS.  This is a similar concept to cutting-plane algorithms, where difficult constraints, relatively few of which are active, are relaxed and individually added as they are violated.  A CT must be used where a \textit{choice} of constraints to add occur when a violation occurs.  In case of MAPF, given a conflict between agents, which agent to be rerouted is not known, and thus branching is used evaluate rerouting each UAV in conflict.  Full details of CBS are given in \cite{sharon2015conflict}.  

\begin{figure}
    \vspace{.2cm}
    \centering
    \includegraphics[width = 0.45\textwidth]{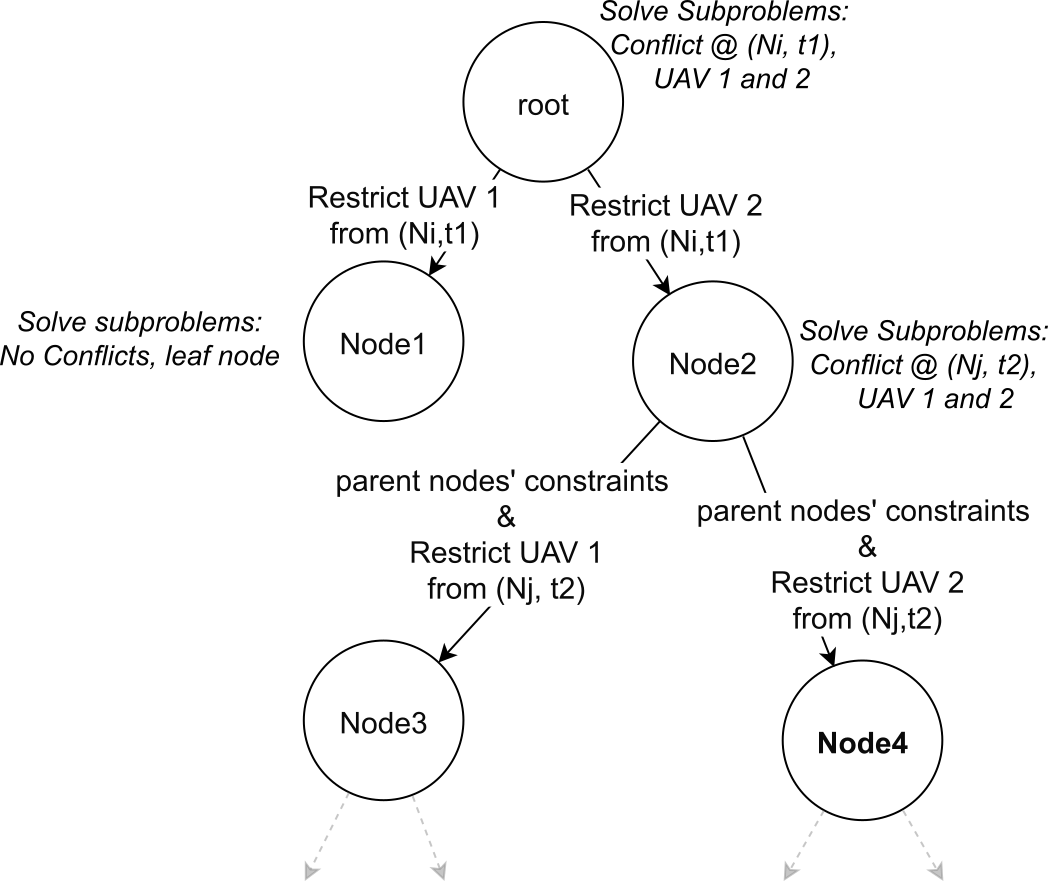}
    \caption{Conflict Based Search - Partial CT Example}
    \label{fig:CBS}
\end{figure}

Fig. \ref{fig:CBS} shows the basic process of CBS, using a CT specific to the NRHF-MAPF.  The root node of the CT consists of solving the NRHF-MAPF with only Equations \eqref{MILPobj}-\eqref{MILPtime2} applied, which is a set of NRHFSPP's, one for each UAV.  All other nodes will contain conflict constraints, which prevent specific UAVs from occupying nodes/edges only at specific times.  At each iteration of CBS, the minimum cost node, according to the subproblem solution cost, is selected.  If the minimum cost node in the CT is feasible, then it is returned at the optimal.  Otherwise, the conflicts occurring in the selection node's subproblem solution are dealt with.  The first conflict in time is picked, and for each UAV in that conflict a child node is produced which restricts all other UAVs from occupying the conflict zone (location/time pair).  Each node contains the conflict constraints relevant to the conflict in the parent node's solution, as well as the constraints from all it's parents up to the root node. 
This process occurs iteratively until the optimal solution is found, or the CT is searched exhaustively.  

As stated, the root node is simply a set of NRHFSPP's, however all other nodes contain conflict constraints which prevent UAVs from occupying a node or edge \textit{at a specific time}.  This is different than the standard NRHFSPP, but can be easily dealt with, discussed in Section \ref{sec:labelalgo}.

\section{Temporal Labeling Algorithm}\label{sec:labelalgo}
To solve the NRHF-MAPF using CBS, the subproblem solved at each node is a set of (NRHFSPP) with the additional conflict constraints. We refer to these conflict constraints as dynamic obstacles, as they function as obstacles which are only present at certain times. These dynamic obstacles must also be dealt with in the SPP subproblems of standard  MAPF.  The dynamic obstacles are handled in the NRHFSPP in a similar manner as how they are dealt with in SPP.  Unlike the SPP, the NRHFSPP is NP-Hard \cite{scott2024power} and does not scale well past a critical problem size.  However, as presented in \cite{scott2024power}, the NRHFSPP can be solved quickly with a labeling algorithm for graphs of tens of thousands of nodes.  For graphs up to some size, we expect the NRHF-MAPF to be tractable despite the NP-Hard subproblem, where the primary bottleneck in most cases is number of agents, which does not directly affect the NRHFSPP tractability.  This is considered empirically in Section \ref{sec:results}.  

We briefly explain the basis of the labeling algorithm here, while it is presented in detail in \cite{scott2024power}.  When solving SPP with Dijkstra's algorithm \cite{dijkstra1959note} or $A^*$ \cite{hart1968formal}, there is only a single value (path cost) associated with each node, which is replaced as better paths are found. However, in the case of resource-constrained SPP and NRHFSPP,  a single value cannot be kept for each node but rather the path cost and all resource values.  Thus, a list of labels is kept for each node, where each label corresponds to a different path to reach that node and thus different resource consumption.  One path may have lower cost while another may have lower resource consumption.  A label is removed from a node's list only if it is dominated by another label, meaning it is worse with respect to both the cost and the resource consumption in comparison to the dominating label. 

Labels are iteratively added to and removed from an open list.  At each iteration of the search, the label of minimum cost is selected from the open list.  The chosen label is expanded by extending its associated partial path to all neighboring nodes of the node which the chosen label is associated with.  These new labels are then added to the open list if they are undominated with respect to the current labels in both the open list and the closed list.  The chosen label is then added to a closed list, which consists only of the labels which have already been treated. The algorithm terminates when the exploration reaches the goal node, based on by a lower-bound cost-to-go heuristic similar to that of $A^*$.

This labeling algorithm to solve the NRHFSPP addresses conflict constraints by treating them as dynamic obstacles.  This is dealt with in standard MAPF by searching SPP over an auxiliary graph.  We follow a similar suite and by modifying the labeling algorithm to solve the NRHFSPP with the dynamic obstacles. This is done by creating an auxiliary graph which relates the connectivity between \textit{states}, defined as (\textit{node, time}) pairs.  Movement between locations also involves movement through time. Thus, occupying the same location at different times is different \textit{states}, and therefore two different nodes in the auxiliary graph.  The path planning search is then posed on this auxiliary graph.  This entails finding a minimum cost path through the auxiliary graph from the start state (with some defined start time) to any state for which the location is the goal node on the original graph. The size of the complete manifestation of this auxiliary can be orders of magnitude larger than the original graph. However, nodes in the auxiliary graph can be created implicitly as this graph is searched.  A state is only added explicitly to the auxiliary graph when a path from a neighboring state is expanded into it.  This implicit generation of the graph keeps the approach tractable, as on average only a small subset of states need to be generated explicitly .  

%%%%%%%%%%%%%%%%%%%%%%%%%%%%%%%%%%%%%%%%%%%%%%%%%%%%%%%%%%%%%%%%%%%%%%%%%%%% 
The pseudo-code of the temporal labeling algorithm is presented in Algorithm 1, where $EFF$ is a function which determines efficiency of a label across the open and closed sets, as discussed in \cite{scott2024power}.  The  minimum cost label is picked in line 4, according to an $f$-cost, which includes an estimated cost-to-go as done in \cite{scott2024power}.  Two different labels are created when exploring the neighboring nodes: one with generator on and one with generator off, shown in lines 8 and 13.  The feasibility checks in lines 11 and 17 address the conflict constraints present in CBS nodes.

\begin{algorithm}
	\SetKwInOut{Input}{input}
	\SetKwInOut{Output}{output}
	
	\Input{Time Cost Matrix: $T \in \R^{N \times N}$ \\ Edge Cost Matrix: $D \in \R^{N \times N}$ \\ Noise Restriction Matrix: $G \in \R^{N \times N}$ \\ Adjacency Matrix: $A \in \R^{N\times N}$\\ Energy Cost Matrix: $C \in \R^{N \times N}$ \\ Energy Transfer Matrix: $Z \in \R^{N \times N}$ \\ Initial Battery/Generator States: $B_0, Q_0$\\Starting node index: $S$\\ Final node index: $F$ \\MAPF constraints for agent:  $M$  \\}
	\Output{Minimum Cost Path and corresponding Generator Pattern for Single Agent}
	
	\SetAlgoLined
	$H  \gets  \{ (d=0,t=0, B_0, Q_0) \}$  // Open Set \\
	% $H_i \gets \emptyset  \quad \forall i \in N, i \neq S$ \\
	$P \gets \emptyset  \quad \forall i \in N$  // Closed Set \\ 
	
	\While{$H \neq \emptyset $}{
		
		% FT = \textproc{Minimum}(H) \\
		FT = $\argmin_{f(\cdot)} H$ \\
		// Label FT: $(d_i^k, t_i^k, b_i^k, q_i^k)$ for node $i$ and associated path $X_k$, $Y_k$ \\
		\For{$(j, t_j' = t_i^k + t_{ij}) \in $  successors$(i, t_i^k$)}{ \label{alg1:nbrnodes}
			// Extend with generator on: \\
			$X_a = [X_k, j]$, $Y_a = [Y_k, 1]$ \\
			$(d_j^a, t_j^a, b_j^a, q_j^a) = (d_i^k  + D_{ij}, t_j', b_i^k - C_{ij} + Z_{ij}, q_i^k - Z_{ij})$ \label{alg1:labelgenon}\\
			\uIf{Feasible($X_a$, $Y_a$) \textbf{and} $(d_j^a, t_j^a) \notin M$ \textbf{and} $j \not\in X_a$ \textbf{and} EFF($(d_j^a, t_j^a, b_j^a, q_j^a), H)$ }{ \label{alg1:feascheck1}
				$H \gets H \cup (d_j^a, t_j^a, b_j^a, q_j^a)$
			}
			// Extend with generator off: \\
			$X_b = [X_k, j]$  \\ 
			$Y_b = [Y_k, 0]$ \\
			$(D_j^b, t_j^b, b_j^b, q_j^b) = (D_i^k  + D_{ij}, t_j', b_i^k - C_{ij}, q_i^k)$ \label{alg1:labelgenoff} \\
			\uIf{$Feasible$($X_b$, $Y_b$) \textbf{and}$(d_j^b, t_j^b) \notin M$ \textbf{and} $j \not\in X_b$ \textbf{and} EFF($(D_j^b, t_j^b, b_j^b, q_j^b), H)$}{ \label{alg1:feascheck2}
				$H \gets \cup (D_j^b, t_j^b, b_j^b, q_j^b)$
			}
		} %end for
		% Remove $FT$ from untreated set $H$\\
		% Add $FT$ to treated set $P$
		$H \gets H \setminus FT$ \\
		$P \gets P \cup FT$
		
	} %end for 
	
	\caption{Single Vehicle Labeling Algorithm} \label{algo:labeling}
	
\end{algorithm}

\section{Results}\label{sec:results}
The CBS algorithm and the labeling algorithm adapted to address dynamic obstacles were implemented in Julia language v1.9 \cite{bezanson2012julia}.  The code for CBS was based on that presented in \cite{choudhury2021efficient}; here, it is extended to deal with the NRHF-MAPF.  All computational experiments were run on a Windows 10 machine with an Intel 4th-generation i5 processor and 16GB of RAM.  We present the results of test problems for a varying number of nodes and number of UAVs.  

Instances were generated by taking the grid problems tested in \cite{scott2024power} and adding random start and goal locations for each agent.  Three problem sizes are tested: 5x5 grid with 4 UAVs, 10x10 grid with 5 UAVs, and 15x15 grid with 10 UAVs. For each size, 50 instances were generated and tested.  Start locations are picked such that no two agents share the same starting location.  All UAVs start their paths at the same time.  The goal location for a UAV is ensured to never be the start location of that same UAV, but it can be same as the start location of other UAVs.  When generating random instances, often the problem is able to be solved at the root node of the CBS.  That is, solving each agent's path individually with no added constraints results in a feasible MAPF solution.  These cases are ignored and not included in our test set, as they are trivial to solve and do not give meaningful insight into performance of our algorithms.  Similarly, it may be the case that a MAPF problem is infeasible, and this is exacerbated by the power and noise constraints present in the problem.  A maximum time limit was set at $120$ seconds, such that the problems that are not solved within this time limit are assumed to be infeasible and are removed from the test problems. At this maximum time limit, tens of thousands of CT nodes have been searched, although this is not a true indication of infeasibility, and more rigorous testing would include increasing this time-limit.

\begin{figure}
    \centering
    \includegraphics[width = 0.4\textwidth]{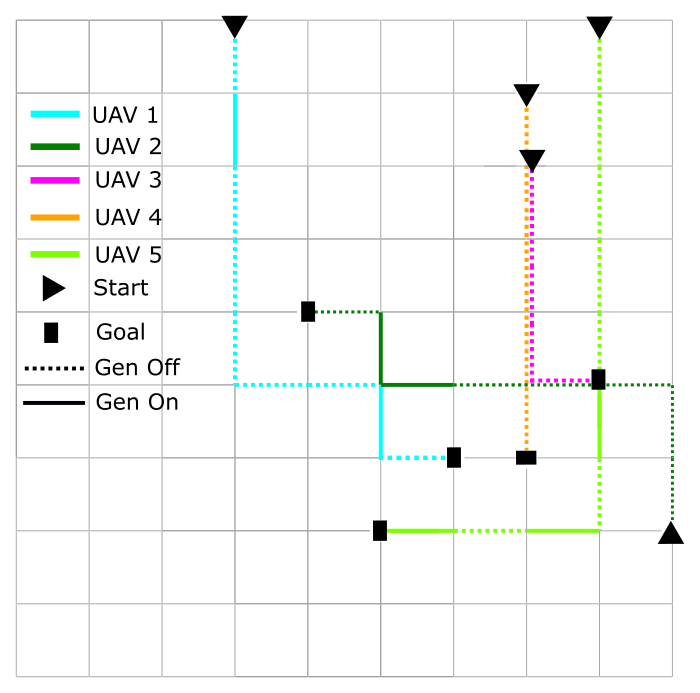}
    \caption{Relaxed Problem (No Conflict Constraints)- 5 UAVs - 100 Nodes - Noise Restricted Zones Hidden - 1 Time Step per Edge - Uniform Start Time}
    \label{fig:paths_init}
\end{figure}

\begin{figure}
    \centering
    \includegraphics[width = 0.4\textwidth]{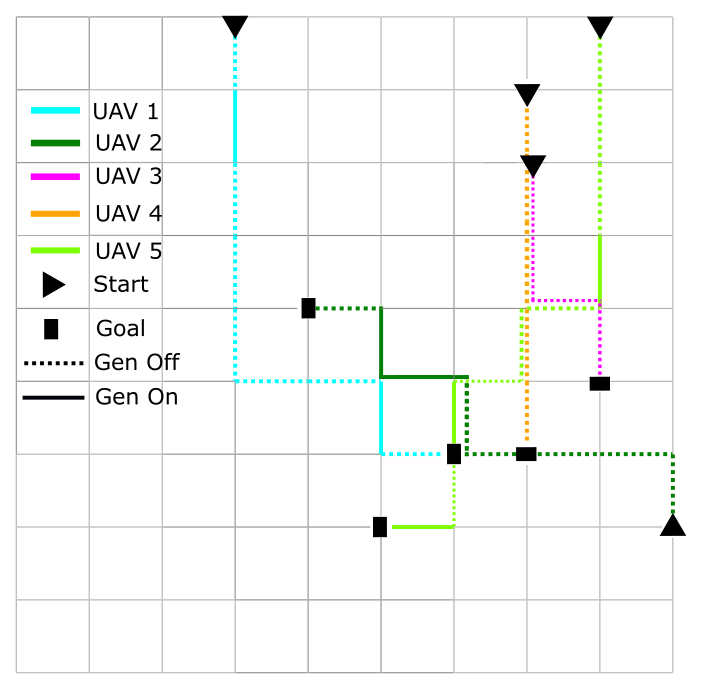}
    \caption{MAPF CBS Solution - 5 UAVs - 100 Nodes - Noise Restricted Zones Hidden - 1 Time Step per Edge - Uniform Start Time}
    \label{fig:paths_final}
\end{figure}

\begin{figure}
    \centering
    \includegraphics[width=.9\columnwidth]{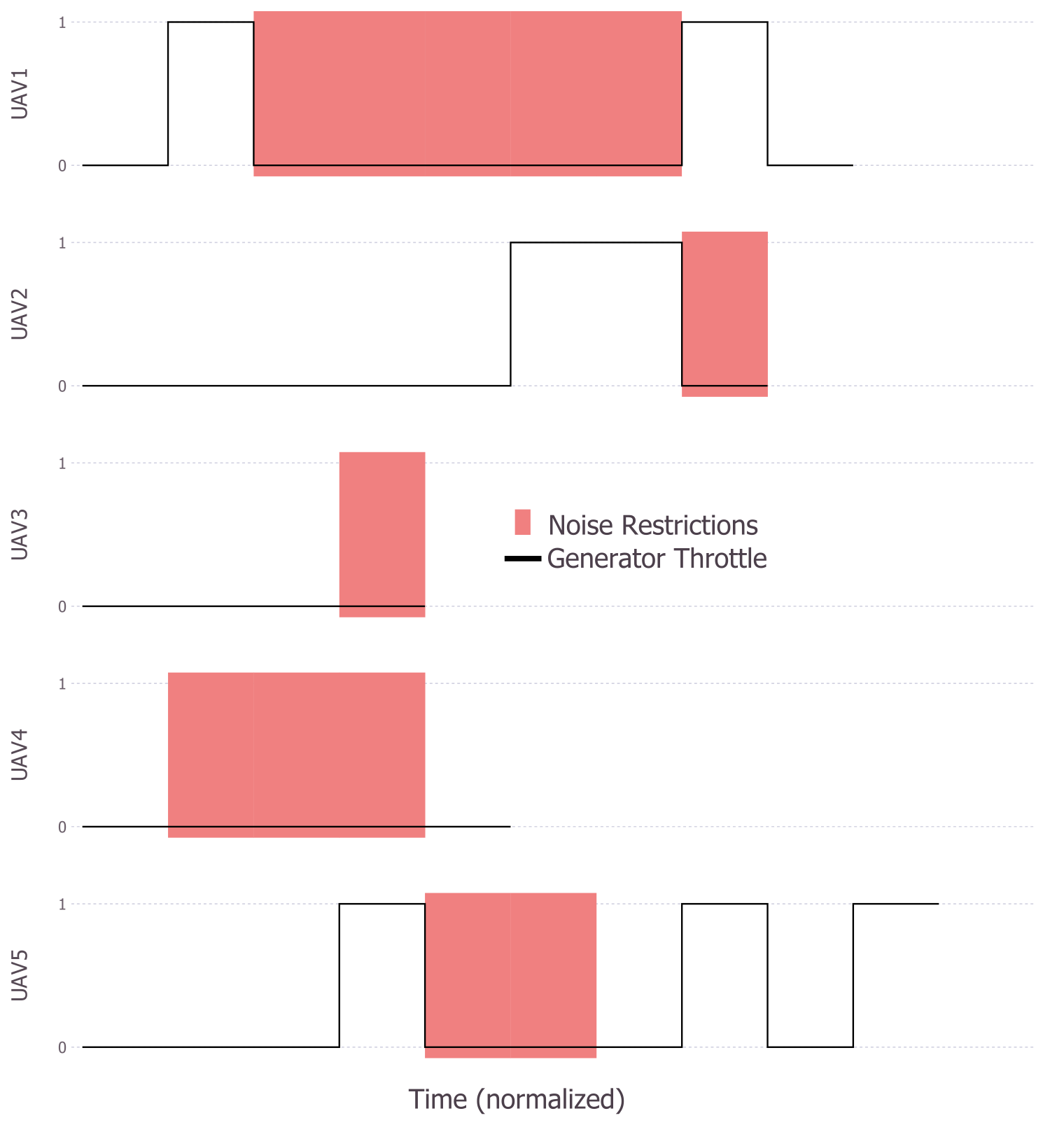}
    \caption{Generator Patterns from CBS Solution - 5 UAVs - 100 Nodes }
    \label{fig:gens}
\end{figure}

For an example problem instance, the solution produced at the root node of the CT, where no conflicts are prevented, is shown in Fig. \ref{fig:paths_init}. One can observe that there are conflicts between many of the paths in the middle-right of the graph.  Here, $37$ CT nodes were explored by the CBS before the optimal solution was found. The optimal solution, with no path conflicts, and the corresponding generator schedule are presented in Figs. \ref{fig:paths_final} and \ref{fig:gens}, respectively. The noise restrictions present along each of the agents' paths are given in Fig. \ref{fig:gens}, shown as red shaded region.


\begin{figure}
    \centering
    \includegraphics[width = 0.45\textwidth]{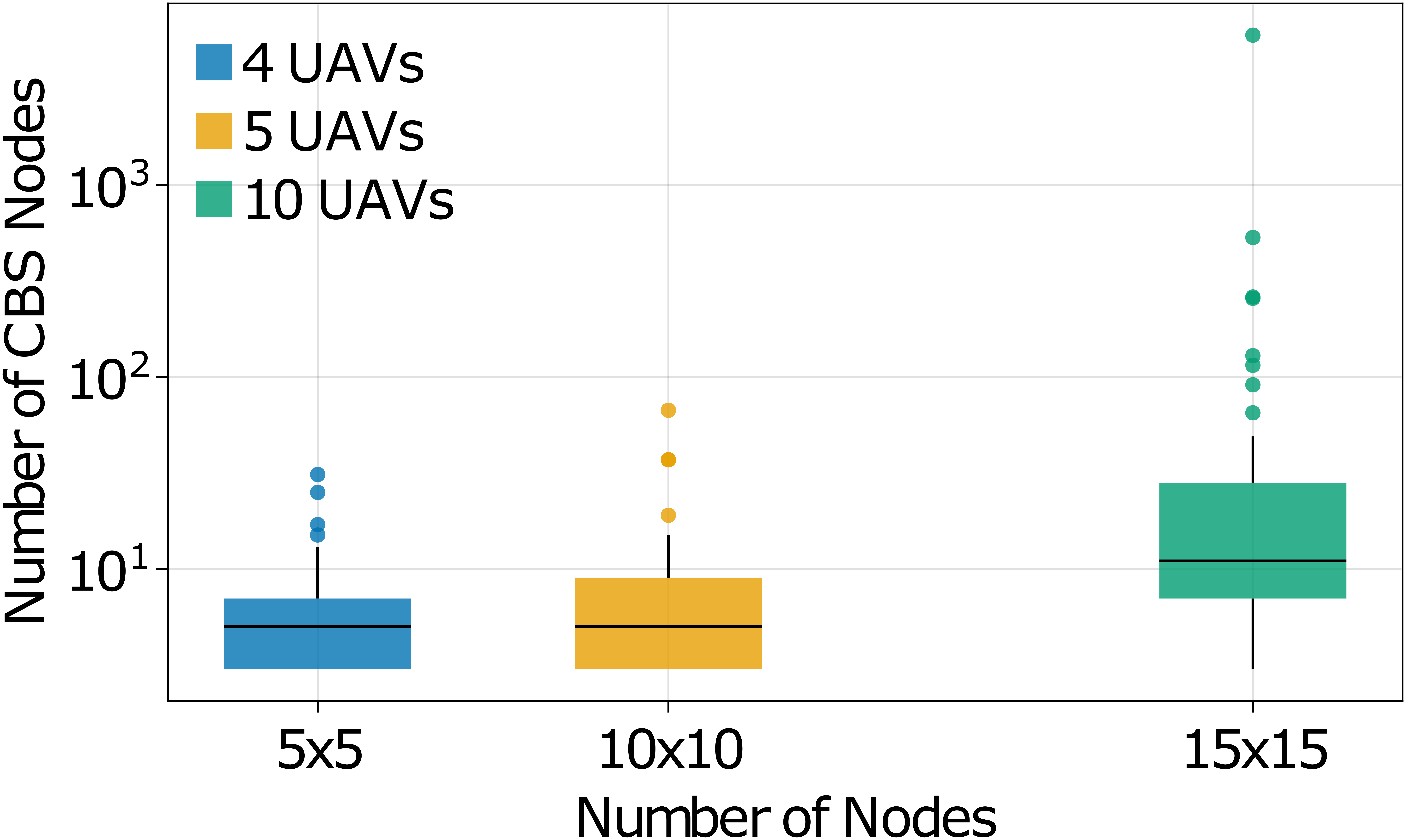}
    \caption{CBS Nodes Searched vs Problem Size}
    \label{fig:node_scaling}
\end{figure}

\begin{figure}
    \vspace{.2cm}
    \centering
    \includegraphics[width = 0.45\textwidth]{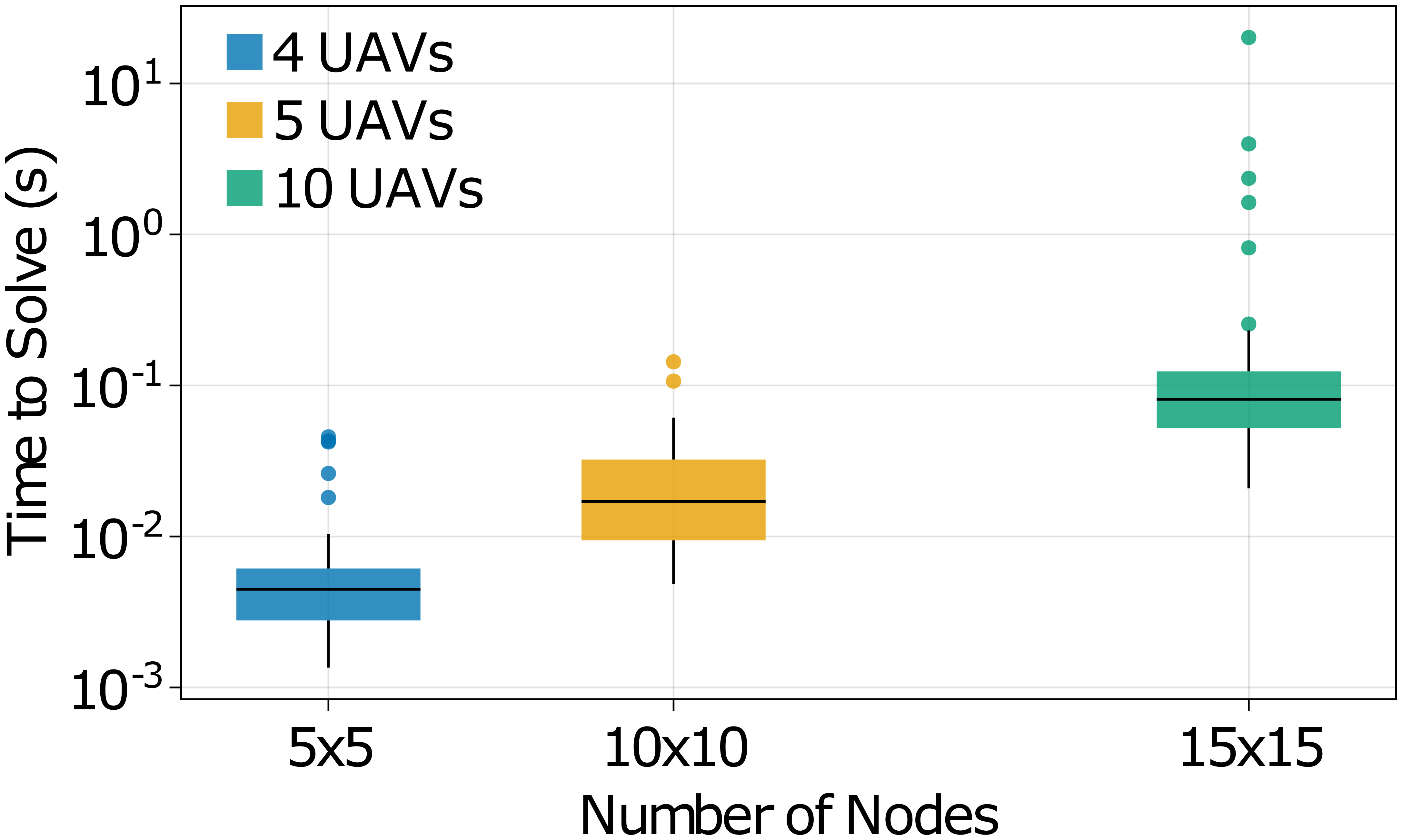}
    \caption{Time to Solve (s) vs Problem Size}
    \label{fig:time_scaling}
\end{figure}

The number of CBS nodes searched in each instance vs the problem size is presented in Fig. \ref{fig:node_scaling}.  The time to solve vs problem size is given in Fig. \ref{fig:time_scaling}.  In each, it can be seen that the difficulty to solve increases with graph size and number of agents.  As seen in the Figures, problems up to 225 nodes and 10 agents can be solved to optimality quickly on average, with the hardest problem taking just over 20 seconds with a total of 6041 CT nodes searched in CBS. 

\begin{figure}
    \centering
    \includegraphics[width = 0.45\textwidth]{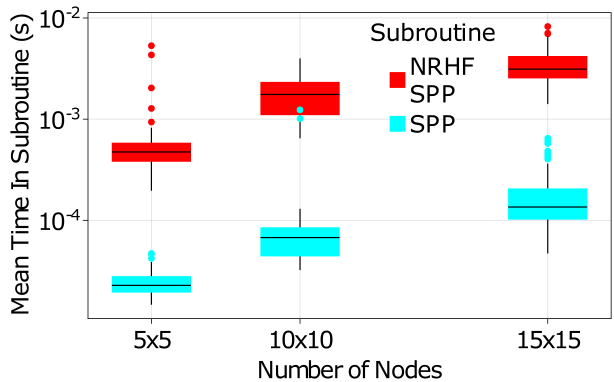}
    \caption{Subproblem Comparison: SPP vs NRHFSPP}
    \label{fig:subroutine_comp}
\end{figure}

As stated above, the subproblem being solved is NP-Hard; however, the labeling algorithm is able to quickly solve for all the instances tested here.  We compare the time to solve the NRHFSPP subproblem with the same single-agent instance with the power and noise constraints removed.  Without these constraints, the problem reduces to a standard SPP, which we solve using $A^*$.  The comparison of the computation times is shown in Fig. \ref{fig:subroutine_comp}.  As expected, the $A^*$ algorithm is faster by orders of magnitude.  It can also be seen in Fig. \ref{fig:subroutine_comp} that there are no especially difficult instances for the NRHFSPP.  This indicates the bottleneck of the CBS search is due to the larger CT rather than the labeling algorithm.  This is reflected in Figures \ref{fig:node_scaling} and \ref{fig:time_scaling}, where the problems which took the longest to solve also had an especially large number of CBS nodes searched. 
Thus it can be said that the especially difficult problems are difficult in terms of finding the optimal \textit{set} of conflict-free paths rather than the individual paths themselves.

\section{Conclusions}\label{sec:conclusions}
Presented was a novel MAPF variant for a fleet of hybrid-fuel UAVs subject to area-dependent noise restrictions. We utilized the CBS algorithm which includes iteratively solving single agent path planning problem with deconflicting constraints added iteratively.  The results show that the NRHF-MAPF can be effectively solved using CBS with the NRHFSPP as a subproblem.  Although the subproblem NRHFSPP is NP-hard, the labeling algorithm presented performs well for the instances tested. 

Extensions of this work include testing the NRHF-MAPF for larger graphs and larger number of agents. Another area of extension includes the implementation of ECBS or other heuristics, which trade optimality guarantees in exchange for improved time-to-solve.  This can also be extended to the single-agent problem, where a heuristic search is developed to the end of quicker solutions for the NRHF-MAPF.

\bibliographystyle{IEEEtran}
\bibliography{bibfile}

% Generated by IEEEtran.bst, version: 1.12 (2007/01/11)
\begin{thebibliography}{10}
\providecommand{\url}[1]{#1}
\csname url@samestyle\endcsname
\providecommand{\newblock}{\relax}
\providecommand{\bibinfo}[2]{#2}
\providecommand{\BIBentrySTDinterwordspacing}{\spaceskip=0pt\relax}
\providecommand{\BIBentryALTinterwordstretchfactor}{4}
\providecommand{\BIBentryALTinterwordspacing}{\spaceskip=\fontdimen2\font plus
\BIBentryALTinterwordstretchfactor\fontdimen3\font minus
  \fontdimen4\font\relax}
\providecommand{\BIBforeignlanguage}[2]{{%
\expandafter\ifx\csname l@#1\endcsname\relax
\typeout{** WARNING: IEEEtran.bst: No hyphenation pattern has been}%
\typeout{** loaded for the language `#1'. Using the pattern for}%
\typeout{** the default language instead.}%
\else
\language=\csname l@#1\endcsname
\fi
#2}}
\providecommand{\BIBdecl}{\relax}
\BIBdecl

\bibitem{yedavalli2019assessment}
P.~Yedavalli and J.~Mooberry, ``{An assessment of public perception of urban
  air mobility (UAM)},'' \emph{Airbus UTM: Defining Future Skies}, 2019.

\bibitem{TOWNSEND2020e05285}
\BIBentryALTinterwordspacing
A.~Townsend, I.~N. Jiya, C.~Martinson, D.~Bessarabov, and R.~Gouws, ``A
  comprehensive review of energy sources for unmanned aerial vehicles, their
  shortfalls and opportunities for improvements,'' \emph{Heliyon}, vol.~6,
  no.~11, p. e05285, 2020. [Online]. Available:
  \url{https://www.sciencedirect.com/science/article/pii/S2405844020321289}
\BIBentrySTDinterwordspacing

\bibitem{scott2024power}
D.~Scott, S.~G. Manyam, I.~E. Weintraub, D.~W. Casbeer, and M.~Kumar, ``Noise
  aware path planning and power management of hybrid fuel uavs,''
  \emph{arXiv:2402.17708}, 2024.

\bibitem{stern2019multi}
R.~Stern, N.~Sturtevant, A.~Felner, S.~Koenig, H.~Ma, T.~Walker, J.~Li,
  D.~Atzmon, L.~Cohen, T.~Kumar \emph{et~al.}, ``Multi-agent pathfinding:
  Definitions, variants, and benchmarks,'' in \emph{Proceedings of the
  International Symposium on Combinatorial Search}, vol.~10, no.~1, 2019, pp.
  151--158.

\bibitem{hopcroft1984complexity}
J.~E. Hopcroft, J.~T. Schwartz, and M.~Sharir, ``On the complexity of motion
  planning for multiple independent objects; pspace-hardness of the"
  warehouseman's problem",'' \emph{The international journal of robotics
  research}, vol.~3, no.~4, pp. 76--88, 1984.

\bibitem{reif1979complexity}
J.~H. Reif, ``Complexity of the mover's problem and generalizations,'' pp.
  421--427, 1979.

\bibitem{surynek2010optimization}
P.~Surynek, ``An optimization variant of multi-robot path planning is
  intractable,'' in \emph{Proceedings of the AAAI conference on artificial
  intelligence}, vol.~24, no.~1, 2010, pp. 1261--1263.

\bibitem{yu2013structure}
J.~Yu and S.~LaValle, ``Structure and intractability of optimal multi-robot
  path planning on graphs,'' in \emph{Proceedings of the AAAI Conference on
  Artificial Intelligence}, vol.~27, no.~1, 2013, pp. 1443--1449.

\bibitem{standley2010finding}
T.~Standley, ``Finding optimal solutions to cooperative pathfinding problems,''
  in \emph{Proceedings of the AAAI Conference on Artificial Intelligence},
  vol.~24, no.~1, 2010, pp. 173--178.

\bibitem{sharon2013increasing}
G.~Sharon, R.~Stern, M.~Goldenberg, and A.~Felner, ``The increasing cost tree
  search for optimal multi-agent pathfinding,'' \emph{Artificial intelligence},
  vol. 195, pp. 470--495, 2013.

\bibitem{sharon2015conflict}
G.~Sharon, R.~Stern, A.~Felner, and N.~R. Sturtevant, ``Conflict-based search
  for optimal multi-agent pathfinding,'' \emph{Artificial Intelligence}, vol.
  219, pp. 40--66, 2015.

\bibitem{barer2014suboptimal}
M.~Barer, G.~Sharon, R.~Stern, and A.~Felner, ``Suboptimal variants of the
  conflict-based search algorithm for the multi-agent pathfinding problem,'' in
  \emph{Proceedings of the International Symposium on Combinatorial Search},
  vol.~5, no.~1, 2014, pp. 19--27.

\bibitem{scott2022power}
D.~Scott, S.~G. Manyam, D.~W. Casbeer, M.~Kumar, M.~J. Rothenberger, and I.~E.
  Weintraub, ``Power management for noise aware path planning of hybrid uavs,''
  in \emph{2022 American Control Conference (ACC)}.\hskip 1em plus 0.5em minus
  0.4em\relax IEEE, 2022, pp. 4280--4285.

\bibitem{manyam2022path}
S.~G. Manyam, D.~W. Casbeer, S.~Darbha, I.~E. Weintraub, and K.~Kalyanam,
  ``Path planning and energy management of hybrid air vehicles for urban air
  mobility,'' \emph{IEEE Robotics and Automation Letters}, vol.~7, no.~4, pp.
  10\,176--10\,183, 2022.

\bibitem{scott2022development}
D.~Scott, S.~G. Manyam, D.~W. Casbeer, M.~Kumar, I.~E. Weintraub, and M.~J.
  Rothenberger, ``Development of linear battery model for path planning with
  mixed integer linear programming: Simulated and experimental validation,''
  \emph{arXiv preprint arXiv:2211.09899}, 2022.

\bibitem{scott2023OCP}
D.~Scott, S.~G. Manyam, D.~W. Casbeer, M.~Kumar, and I.~E. Weintraub, ``Optimal
  generator policy for hybrid fuel uav under airspace noise restrictions,''
  \emph{to be published in MECC 2023}, 2023.

\bibitem{jadischke2023optimal}
J.~H. Jadischke, M.~Wolff, J.~Zumberge, B.~Hencey, and A.~Ngo, ``Optimal route
  planning and power management for hybrid uav using a* algorithm,'' in
  \emph{AIAA AVIATION 2023 Forum}, 2023, p. 4508.

\bibitem{williams2013model}
H.~P. Williams, \emph{Model building in mathematical programming}.\hskip 1em
  plus 0.5em minus 0.4em\relax New York: John Wiley \& Sons, 2013.

\bibitem{dijkstra1959note}
E.~W. Dijkstra \emph{et~al.}, ``A note on two problems in connexion with
  graphs,'' \emph{Numerische mathematik}, vol.~1, no.~1, pp. 269--271, 1959.

\bibitem{hart1968formal}
P.~E. Hart, N.~J. Nilsson, and B.~Raphael, ``A formal basis for the heuristic
  determination of minimum cost paths,'' \emph{IEEE transactions on Systems
  Science and Cybernetics}, vol.~4, no.~2, pp. 100--107, 1968.

\bibitem{bezanson2012julia}
J.~Bezanson, S.~Karpinski, V.~B. Shah, and A.~Edelman, ``Julia: A fast dynamic
  language for technical computing,'' \emph{arXiv preprint arXiv:1209.5145},
  2012.

\bibitem{choudhury2021efficient}
S.~Choudhury, K.~Solovey, M.~J. Kochenderfer, and M.~Pavone, ``Efficient
  large-scale multi-drone delivery using transit networks,'' \emph{Journal of
  Artificial Intelligence Research}, vol.~70, pp. 757--788, 2021.

\end{thebibliography}

\end{document}